\documentstyle[12pt,amssymb]{amsart}
\newcommand{\stopproof}{\hfill \nobreak\medskip $\blacksquare$ \\
\hspace*{\fill}}
\newcommand{\dom}{\mbox{\rm dom}}

\newcommand{\AND}{\mbox{ \rm and }}

\newcommand{\proof}{{\bf Proof:} \ }

\newcommand{\AAA}{{\Bbb A}}

\newcommand{\BB}{{\Bbb B}}

\newcommand{\restricts}{\! \upharpoonright \!}

\newtheorem{theor}{Theorem}[section]
\newtheorem{propo}{Proposition}[section]
\newtheorem{defin}{Definition}[section]
\newtheorem{corol}{Corollary}[section]
\newtheorem{lemma}{Lemma}[section]
\newtheorem{quest}{Question}[section]

 \title{Reaping Numbers of Boolean Algebras}
\author{Alan Dow} \author{Juris Stepr\={a}ns} \author{Stephen Watson}
\address{Department of Mathematics, York University \\
4700 Keele Street \\ North York, Ontario \\ Canada \ \ \ \ M3J 1P3}
\thanks{The research of all three authors is partially supported by 
NSERC}
\begin{document}
\maketitle
\begin{abstract}
A subset ${\cal A}$ of a Boolean algebra $\BB$ is said to be $(n,m)$-reaped
if there is a partition of unity ${\cal P} \subset {\BB}$ of size $n$
such that $|\{b \in {\cal P}: b \wedge a \neq
\emptyset\}| \geq m$ for all $a\in\cal A$.
The reaping number $\frak{r}_{n,m}(\BB)$ of a Boolean algebra
$\BB$ is the minimum cardinality of a set ${\cal A} \subset {\BB}\setminus 
\{0\}$ such
which cannot be $(n,m)$-reaped.
It is shown that, for each $n \in \omega$, there is a Boolean algebra
 $\BB$ such that
$\frak{r}_{n+1,2}(\BB) \neq \frak{r}_{n,2}(\BB)$.
Also, $\{\frak{r}_{n,m}(\BB) : \{n,m\}\subseteq\omega\}$ consists of 
at most two 
consecutive integers. The existence of a Boolean algebra $\BB$ such that
 $\frak{r}_{n,m}(\BB) \neq \frak{r}_{n',m'}(\BB)$ is equivalent
to a statement in finite combinatorics which is also discussed.

\noindent {\bf Keywords:} Boolean algebra, reaping number, tree, colouring
\end{abstract}
\bibliographystyle{amsplain}
\section{Introduction}
A subset $A$ of a Boolean algebra $\Bbb B$, is said to be {\em
reaped} by $b \in \Bbb{B}$, if, for all $a \in A$, both $a\wedge b$
and $a - b$ are non-zero. For any Boolean algebra $\Bbb B$, the
cardinal invariant $\frak{r}(\Bbb B)$ is defined to be the least
cardinality of a subset of $\Bbb B\setminus \{0\}$ which can not be reaped.
The more familiar cardinal invariant of the
continuum known as the {\em reaping number}, $\frak r$,
is therefore nothing more than
$\frak{r}({\cal P}(\omega)/[\omega]^{<\aleph_0})$.
This cardinal invariant was introduced by A. Beslagi\'{c} and E. van Douwen
\cite{BvD99} who showed that the equality $\frak{r} = {\frak c}$ implies that
$\omega^* \setminus \{p\}$ is non-normal for
each ultrafilter $p \in \omega^*$.
The reaping number appears again in a paper of B. Balcar and P. Simon
\cite{ba.si.1} concerning the
$\pi$-character of points in $\omega^*$, the Stone space of
${\cal P}(\omega)/[\omega]^{<\aleph_0}$.
The $\pi$-character of a point $x$ in a topological space $(X,\tau)$ is the
least 
cardinality of a set $\cal{A}\subset \tau$ such that for
every  neighbourhood $V$ of $x$ there is $A\in \cal{A}$ such that $A\subseteq 
V$
--- it is not required that ${\cal A}$ consist of neighbourhoods of $x$.
In \cite{ba.si.1} $\frak r$ is  characterised
 as the minimum of all $\pi$-characters
of points in $\omega^*$. This leads naturally to the question of
whether or not $\frak{r}(\BB)$ is the minimum of all
$\pi$-characters of ultrafilters on $\BB$ when these are considered as points
in the Stone space of $\BB$.

Given a Boolean algebra $\Bbb B$ and an integer $n$, an $n$-partition
of $\Bbb B$ is a set $R \in  [\Bbb{B}]^n$ such that $\bigvee R = 1$ and
$a \wedge b = 0$ for $\{a, b\} \in [R]^2$.  A set $A \subseteq \Bbb B$ will
be said to be $(n,k)$-reaped by $R$ if $R$ is an $n$-partition of $\Bbb B$
such that $|\{ b \in R : b \wedge a \neq 0\}| \geq k$ for all $a \in A$.
Therefore a set $A$ is reaped by $b$ if and only if it is $(2,2)$-reaped
by $\{ b,1-b\}$.  The cardinal $\frak{r}_{n,k}(\Bbb{B})$ is the least
cardinality of a subset of $\Bbb{B}\setminus\{0\}$ 
which cannot be $(n,k)$-reaped
by some $n$-partition of $\Bbb{B}$. Define $\frak{r}_n(\BB)$ to be
$\frak{r}_{n,2}(\BB)$ and the phrase $n$-reaped will be used to mean
$(n,2)$-reaped. Let $\frak{r}_\omega(\Bbb B)$ denote the least cardinal of
a set $A\subseteq \BB\setminus\{0\}$ 
which cannot be $n$-reaped for any $n$.
The sequence of inequalities
\[\frak{r}_2 \leq \frak{r}_3 \leq \frak{r}_{4} \leq \cdots \leq \frak{r}_\omega\]
is easily verified. With these definitions, the results  of
\cite{ba.si.1} can now be stated.
The first establishes a connection with topology.
\begin{propo}
For every  Boolean algebra $\Bbb B$, $\frak{r}_\omega(\Bbb B)$ is equal to the
\label{balc:1}
least cardinal $\kappa$ such that there is an ultrafilter on $\Bbb B$
with $\pi$-character $\kappa$.
\end{propo}
{\bf Proof: }
Note that a $\pi$-base for an ultrafilter cannot be
$n$-reaped for any $n \in \omega$ so it suffices to show that there
is some ultrafilter on $\Bbb B$ whose $\pi$-character is
$\frak{r}_{\omega}(\Bbb B)$.
Let $A \subseteq \Bbb B$ be a set which cannot be $n$-reaped for
any $n\in\omega$ and such that $|A| = \frak{r}_{\omega}(\Bbb B)$.
Then  $G = \{b \in \Bbb{B} : (\forall a \in A) (a- b \neq 0)\}$
generates a proper ideal because otherwise
there is, for some $n\in\omega$, an $n$-partition  consisting of
 elements of $G$ and this
contradicts that $A$ is not  $n$-reaped. This ideal can be extended to a
maximal one and its dual filter is the desired filter. It can be checked
that  every member of this dual filter contains some member of $A$.
\stopproof

While Proposition \ref{balc:1} illustrates the importance of the cardinal
invariants $\frak{r}_n(\Bbb B)$, the next result from \cite{ba.si.1}
 indicates that, in some
cases, they are all equal.
It also shows that the equivalence of $\frak{r}$ and the least $\pi$-character
of a point in $\omega^*$ is due to the incidental fact that
$\cal{P}(\omega)/[\omega]^{<\aleph_0}$ is homogeneous.
\begin{propo}
If $\Bbb B$ is \label{balc:2}
homogenous or complete, then $\frak{r}_2(\Bbb{B}) = \frak{r}_\omega(\Bbb B)$.
\end{propo}
{\bf Proof: }
First suppose that $\Bbb B$ is homogeneous. In fact, simply assume that for
all $b \in \Bbb B$, $\frak{r}_2(\Bbb B) = \frak{r}_2 (\Bbb{B}\restricts b)$.
 Proceed by induction on $n$ to show that if $\BB$ is any Boolean algebra
satisfying these assumptions then $\frak{r}_2(\BB) =\frak{r}_n(\BB)$ ---
this clearly suffices. The case $n=2$ is trivial.
Let $A\subseteq\BB$ be a set of cardinality $\frak{r}_2(\BB)$ which can not
be 2-reaped. 
 It may, without loss of generality be assumed that 
$A$ is closed under
finite products since  closing it off under this operation will
not change its cardinality.
 If $R $ is an $n+1$-partition of
$\Bbb B$ which $n+1$-reaps $A$ fix any
$r \in R$. Since $A$ is not $2$-reaped, there is an $a \in A$,
such that $a \leq r$  or $a \wedge r = 0$.
The first case contradicts that $A$ is $n+1$-reaped by  $R$.
Therefore $A \cap \Bbb{B}\restricts a$ is $n$-reaped  by
$ R \setminus \{r\}$. This contradicts the  induction hypothesis and the 
homogeneity 
of $\BB$.

Now suppose that $\Bbb B$ is complete and choose $A \subseteq \Bbb B$ which
cannot be $2$-reaped and is
of minimal cardinality.  Find a set  $C \subseteq \Bbb B$
such that  
\begin{itemize}
\item $\bigvee C = 1$ \item $c \wedge c' = 0 $ for
all $c$ and $ c'$ in $C$ provided that $c\neq c'$
\item if   $c \in C$ and $b \leq c$ then $\frak{r}_2(\Bbb{B}\restricts b)
= \frak{r}_2(\Bbb{B} \restricts c) $
\end{itemize}
 If $\frak{r}_2(\Bbb{B}\restricts c) = \frak{r}_2(\Bbb{B})$
 for some
$c \in C$ then the homogeneity assumptions of the first part are satisfied.
 Otherwise, for each $c \in C$, choose
$b_c < c$ such that  $A\cap
\Bbb{B}\restricts c$ is 2-reaped by $\{b_c, c-b_c\}$. 
It is easily seen that the partition $\{\bigvee_{c\in C}b_c,
1 - \bigvee_{c\in C}b_c\}$ also $2$-reaps $A$.\stopproof

The question of  whether
$\frak{r}_{2}({\Bbb B}) = \frak{r}_{\omega}({\Bbb B})$
for all Boolean algebras was raised by the authors of \cite{ba.si.1}. 
It will be shown that the answer is no.
The same result was obtained independently by B. Balcar and P. Simon
in \cite{ba.si.2}. The next section contains results which 
yield this as a corollary and establish the best possible upper bound for $\frak{r}_n(\BB)$ for $n\leq \omega$. The more general cardinal invariant $\frak{r}_{i,j}(\BB)$ which 
is introduced in this paper plays a central role in this analysis.
\section{Reaping Numbers
and a Finite Combinatorial Property}
This section is devoted to establishing an equivalence between an 
assertion abount trees and 
the statement that $\frak{r}_{i,j}(\BB)\leq \frak{r}_{k,n}(\BB)$ for
all Boolean algebras $\BB$.
Throughout this paper the term {\em tree} will be reserved for 
sets of sequences of integers which are closed under initial segments. 
The canonical examples of trees
are the sets $T(D,m) = \cup_{\xi\in D}{}^{\xi+1\cap D}m$ 
where $D$ is a set of ordinals 
--- in other words, $T(D,m)$ is the set of all $m$-valued functions with domain
an initial segment of $D$. 
All trees considered in this paper will be subsets of some $T(D,m)$. 
The ordering
on a tree is the restriction to an initial segment relation; in other words,
if $\sigma$ and $\theta$ belong to $T(D,m)$ then $\sigma \leq \theta$ if
and only if there is some $\delta\in D$ such that $\sigma = \theta
\restriction \delta$. The set of maximal elements of a tree $T$ will be denoted by
$\mu(T)$. A tree $T$ is said to be $k$-branching if every
$\sigma \in T\setminus \mu(T)$ has exactly $k$ successors.
\begin{defin}
The quadruple of integers $(i,j,k,m)$ will be said to satisfy
the  property $P(i,j,k,m)$
if and only if for every finite, $i$-branching
tree $T$ and every colouring $\chi:\mu(T)\rightarrow j$
 of the maximal members of $T$
by $j$ colours, there is a $k$-branching subtree $T'\subseteq T$ such
that  $\mu(T')\subseteq \mu(T)$ and  the range of $\chi\restriction\mu(T')$
has no more than $m$ colours.
\end{defin}
The main result of this section is Theorem~\ref{ijkl}
which establishes  a relationship between the assertion 
that $\frak{r}_{i,j}(\BB)\leq \frak{r}_{k,m}(\BB)$ for every Boolean
algebra $\BB$
and the property $P(i,k,j,m)$. The next lemma is the key to establishing this theorem.
\begin{lemma}
For all pairs of integers $i$ and $k$
such that $k\leq i$ there \label{keyijkl}
 are $k$-branching trees $\{T_{n}
: n\in\omega\}$ such that
for each $F \in [\omega_{1}]^{< \aleph_0}$ and 
for each  $k$-branching subtree $T\subseteq T(F,i)$ there 
are infinitely many $n \in \omega$ 
such that $T =\{\psi\restriction F :\psi\in
 T_{n}\}$.
\end{lemma}
{\bf Proof: } It is possible to prove this
result by constructing the natural partial order to force the trees $T_n$
with finite conditions and then appeal to D. Velleman's Martin's
Axiom type equivalence to
$\omega$-morasses \cite{vell}. For the reader's convenience however,
a self contained proof will be presented.
The trees $T_n$ will be chosen to form a dense subset of a certain product
space whose points can be interpreted as trees. The fact that the product
of no more than $2^{\aleph_0}$ separable spaces is separable will play an
important role \cite{enge.gt}.

For each $\alpha\in\omega_1$ let $S(\alpha)$ be the set of all
colourings $\chi :\mu(T(F_{\chi},k))\rightarrow [i]^k$ where $F_{\chi}\in
[\alpha]^{<\aleph_0}$ with the discrete topology and let
$S=\prod_{\alpha\in\omega_1}S(\alpha)$. Since each $S(\alpha)$ is
countable it follows that $S$ is separable.
Moreover each element of $S$ can be considered to be a $k$-branching subtree
of $T(\omega_1,i)$. In order to do this let $a^e : a\rightarrow k$ be
the unique, order preserving bijection from $a$ to $k$ for each
$a\in [i]^k$. Then, if $x\in S$  define a $k$-branching subtree
$x^*_{\theta}\subseteq T(\theta,i)$ by induction on $\theta$. If
$\theta$ is limit of countable cofinality then \[x^*_{\theta}
= \{\sigma : (\forall \zeta\in\theta)(f\restriction \zeta \in x^*_{\zeta})\}\]
and $x^* $ will be defined to be $\cup_{\theta\in\omega_1}x^*_{\theta}$.
If $\theta = \zeta +1$ then let $C_{\zeta}:x^*_{\zeta}\rightarrow
T(\zeta,k)$ be the canonical collapse of $x^*_{\zeta}$; namely, letting
$a(\sigma)$ be the set of successors of $\sigma$ in $x^*_{\zeta}$,
$C_{\zeta}(\sigma)(\beta) = a(\sigma\restriction\beta)^e(\sigma(\beta))$.
Then $x^*{\theta}$ is defined to 
be the set of all $\sigma\in T(\theta,i)$ such that
$\sigma\restriction\zeta\in x^*_{\zeta}$ and
$\sigma(\zeta) \in x(\zeta)(C_{\zeta}(\sigma\restriction\zeta)
\restriction F_{x(\xi)})$. 

All that has to be checked is that if $F\in
[\omega_1]^{<\aleph_0}$ and
 $T\subseteq T(F,i)$ is a $k$-branching subtree then
\[\{x\in S : (\forall \sigma\in x^*)(\sigma\restriction F\in T)\}\]
is a nonempty open set. This is a routine exercise.
\stopproof

\begin{theor}
The property \label{ijkl} $P(i,j,k,m-1)$ is satisfied if and only if 
there is a Boolean algebra $\BB$ such that
$\frak{r}_{i,k}(\BB) > \frak{r}_{j,m}(\BB)$.
\end{theor}
{\bf Proof: }
If $P(i,j,k,m-1)$ fails, then there is some integer $K$ and a colouring
$\chi : \mu(T(K,i))\rightarrow j$ 
such that every $k$ branching subtree is coloured by at
 least $m$ colours.
Let $\BB$ be a Boolean algebra and  $F\subseteq \BB$ be
of cardinality less than $\frak{r}_{i,k}(\BB)$. It must be shown that
$F$ can be $(j,m)$-reaped, so suppose that this is not possible.
 
For each $a\in \BB$ it is possible to 
choose an $i$-partition $\{p_s(a) : s\in i\}$
which $(i,k)$-reaps $\{f\cap a : f\in F\AND f\cap a\neq 0\}$. 
For each $\sigma\in T(K,i)$ it is possible to define inductively 
a function $P:T(K,i)\rightarrow\BB$ as follows:
$P(\emptyset) = 1_{\BB}$ and, if the domain of $\sigma$ is $L$, then
 $$P(\sigma\cup\{(L,s)\}) = p_s(P(\sigma))\cap P(\sigma)$$
Define $b_c = \cup\{P(\sigma) :
\sigma :K\rightarrow i\AND \chi(\sigma) = c\}$ for $c\in j$.
 It must be that the partition $\{b_c : c\in j\}$ does not 
$(j,m)$-reap $F$ and so there is some $f$, an element  of $F$,
 which witnesses this fact.
Using the fact that $\{p_s(a) : s\in i\}$ is a partition
which $(i,k)$-reaps $F$ it follows that there is a $k$ branching subtree
$T'\subseteq T(K,i)$ such that
$f\cap P(\sigma)\neq 0$ for each  $\sigma \in \mu(T')$.
 But $\mu(T')$  is coloured by $\chi$ with at least 
$m$ colours
and  $f\cap b_c\neq 0$ if some member of $\mu(T')$ is coloured with colour
$c$. It follows that $\{c\in j : b_c\cap f\neq 0\}$ has at most $m$ elements 
and this
is a contradiction to the choice of $f$.

Suppose $P(i,j,k,m)$ holds. It will be shown that
there is a Boolean algebra in which a countable family cannot
be $(j,m+1)$ reaped but every countable family can be $(i,k)$ reaped.
Let $\Bbb{B} = \Bbb{C}_{\omega}\oplus\Bbb{F}_{\omega_1}$ where
$\Bbb{C}_{\omega}$ is the Boolean subalgebra of $\cal{P}(\omega)$
generated by the finite sets and $\Bbb{F}_{\omega_1}$ is the free
algebra of size $\aleph_1$ represented as $\prod_{\xi\in\omega_1}{\cal P}(i)$.
The algebra $\BB$ is generated by rectangles of the form
$F \times [\sigma]$ where $F$ is a finite or cofinite subset
of $\omega$ and $\sigma$ is a finite partial function from
$\omega_{1}$ into $i$;
these rectangles will be denoted by $F \otimes [\sigma]$.
The trees
$T_{n}$ of Lemma~\ref{keyijkl} will be used
 to generate a quotient algebra
${\Bbb A}$. Let $\cal{G}_n = \{\sigma\restriction F : \sigma \in T_n\AND
 F\in [\omega_1]^{<\aleph_0}\}$. Let $\cal{I}$ be the ideal in $\BB$
 generated $\{\{n\} \otimes [\psi] : \psi \not\in {\cal G}_{n}\}$
and define $\AAA = \BB/\cal{I}$.

To see that $\frak{r}_{i,k}(\Bbb{A}) = \aleph_1$ let
$\Bbb{B}_{\alpha}$ be the subalgebra of $\BB$
 generated by sets of the
form $\{n\} \otimes [\psi]$ where $\dom(\psi)\subseteq \alpha$
and define $\AAA_{\alpha}
 =\Bbb{B}_{\alpha}/\cal{I}$.
It suffices to show that $\AAA_{\alpha}$ is $(i,k)$-reaped by
$\{\omega \otimes [(\alpha,\ell)] : \ell\in i\}$.
Now, if the equivalence class of 
$\{n\} \otimes [\psi]$ is not empty
in $\AAA_{\alpha}$
 then $\psi \in {\cal G}_{n}$.
By Lemma~\ref{keyijkl}, given this $\psi$ and $n$ and $\alpha$,
it follows that there is $\sigma\in T_n$ such that $\psi\subset\sigma$
and there is $a \in [i]^k$
such that $\{\sigma\cup\{(\alpha,t)\}: t \in a\} \subseteq T_{n}$. Hence
  $\{\psi\cup\{(\alpha,t)\}: t \in a\} \subseteq \cal{G}_{n}$. 
>From this it will follow
that $\AAA_{\alpha}$ is $(i,k)$-reaped provided that it can be shown that
every nonempty member of $\AAA_{\alpha}$ contains something of the form
$\{n\}\otimes [\psi]$ such that $\psi\in\cal{G}_n$.

To this end notice that every nonempty member of $\AAA_{\alpha}$ contains
 something of the form
$C\otimes [\psi]$. If $C$ is finite then 
$C\otimes [\psi] = \cup_{n\in C}\{n\}\otimes[\psi]$ and the fact that
$C\otimes [\psi]\neq 0$ implies that there is some $n\in C$ such that
$\{n\}\otimes [\psi]\notin \cal{I}$. If $C$ is infinite then it contain
$\omega\setminus J$ for some $J\in\omega$. Let $\Gamma$ be the domain of $\psi$
and let $T$ be any $k$-branching subtree of $T(\Gamma,i)$ such that $\psi
\in \mu(T)$. From the properties of $\{T_n : n\in\omega\}$ stated in 
Lemma~\ref{keyijkl} it follows that there is some $j > J$ such that
$\psi\in T\subset \cal{G}_j$. Hence $\{j\}\otimes[\psi]\subseteq
C\otimes[\psi]$.

To see that $\frak{r}_{j,m+1}(\Bbb{A}) = \aleph_0$,
it will be  shown that $\{n\} \otimes [\emptyset]$ is a countable family
which cannot be $(j,m+1)$ reaped.
Note that each $a\in {\Bbb B}$ which is not zero in
$\Bbb A$ can be represented as
\[a = \bigvee_{\psi\in T(F,i)} C(\psi) \otimes [\psi]\]
where $F \in [\omega_{1}]^{<\aleph_0}$ and where
$C(\psi)$ is a finite or cofinite subset of $\omega$.
It can therefore be deduced that any partition of unity
$\{a_z : z \in j\}$ can be represented as
\[a_z = \bigvee_{\psi\in T(F,i)} C_z(\psi) \otimes [\psi]\]
for some $F \in [\omega_{1}]^{<\aleph_0}$
where $\{C_z(\psi) : z \in j\}$ is a partition of $\omega$
into finite and cofinite sets, for each $\psi$.
Given such a partition of unity it is possible to choose an integer
$K$ such that for every $z\in j$ and $\psi \in T(F,i)$ either
$C_z(\psi)\subseteq K$ or $\omega\setminus K\subseteq C_z(\psi)$.

Now define a colouring $\chi$ of $\mu(T(F,I))$ by the rule
$\chi(\psi) = z$ if and only if $\omega\setminus K \subseteq C_z(\psi)$.
Now, using $P(i,j,k,m)$, there is a $k$ branching subtree $\Gamma$ of
$T(F,i)$ such that the range of
$\chi\restriction \mu(\Gamma)$ has size at most $m$.
By the key property of the family of trees stated in Lemma~\ref{keyijkl}
it follows that
there is an integer $n > K$ such that
$\{\sigma\restriction F : \sigma \in T_{n}\} = \mu(\Gamma)$.

Since $n > K$ it follows that $n \in C_z(\psi)$ if and only if
$\chi(\psi) = z$. Also, if $\{n\}\otimes [\psi]\notin\cal{I}$ and
the domain
of $\psi$ is $F$ then
$\psi\in\mu(\Gamma)$. Hence, if $a_z\wedge \{n\}\otimes [\emptyset]\neq
\emptyset$ then
$\chi(\psi) = z$ for some $\psi \in\mu(\Gamma)$.
 Since the range $\chi\restriction\mu(\Gamma)$ has at most $m$ elements
it follows that $z$ can take on at most $m$ values and the proof is complete.
\stopproof

Note that the Boolean algebra constructed in this section has
the countable chain condition. Note also that it is possible to
generalise the construction
for any cardinal $\kappa$ to obtain a Boolean algebra
$\Bbb A$ such that $\frak{r}_{i,k}(\Bbb{A}) = \kappa^{+}$ while
$\frak{r}_{j,m+1}(\Bbb{A}) = \kappa$.
\section{Some Inequalities}
Certain monotonicity results for the cardinal
invariants $\frak{r}_{n,m}(\Bbb{B})$ are easily established.
\begin{lemma}
If \label{mono} $\Bbb{B}$ is any Boolean algebra and $n$ and $m$
are integers with $m \geq n \geq 2$, then the following inequalities hold
\begin{enumerate}
\item $\frak{r}_{m,n}(\Bbb{B}) \geq \frak{r}_{m,n+1}(\Bbb{B})$
\item $\frak{r}_{m+1,n}(\Bbb{B}) \geq \frak{r}_{m,n}(\Bbb{B})$
\item $\frak{r}_{m+1,n+1}(\Bbb{B}) \leq \frak{r}_{m,n}(\Bbb{B})$ and,
more generally,
$\frak{r}_{i,j}(\Bbb{B}) \leq \frak{r}_{m,n}(\Bbb{B})$ provided that
there is some integer $k$ such that $m\leq i\leq mk$ and $j > (n-1)k$
\item $\frak{r}_{m,n}(\Bbb{B}) \geq \frak{r}_{2}(\Bbb{B})$
for all integers $m$ and $n$
\item $\frak{r}_{n,m}(\Bbb{B}) =  \frak{r}_{2}(\Bbb{B})$ for all integers $n$
and $m$  provided that
$n/(m-1)\geq 2$.
\end{enumerate}
\end{lemma}
{\bf Proof: }
The first two assertions are immediate.
To verify the last clause of (3) suppose that
 $\{a_\ell : \ell \in i\}$ is an $i$-partition which
$(i,j)$-reaps $A \subseteq {\Bbb B}$.
Then it is possible to find $\{u_\ell : \ell\in m\}
\subseteq [i]^{\leq k}$ which partition
 $i$ into nonempty sets.  Let $b_{\ell} = \vee\{a_x : x\in
u_{\ell}\}$. It follows that $\{b_{\ell} : \ell \in m\}$ $(m,n)$-reaps
$A$ because of the fact that  $j > (n-1)k$ and so, if less than
$n$ members of $\{b_{\ell} : \ell \in m\}$ meet some member of $A$
then less than $j$ members of $\{a_{\ell} : \ell \in i\}$ meet that
same element.
 
To prove (4) suppose that $A\subseteq \Bbb{B}$ can not be $(i,j)$-reaped
and $|A| < \frak{r}_2(\BB)$. It follows that there is a 2-partition
$\{a_0, a_1\}$ which 2-reaps $A$. Moreover $$|\{a\wedge a_0 : a\in A\}|
\leq |A| < \frak{r}_2(B)$$ and
$|\{a\wedge a_1 : a\in A\}| \leq |A| < \frak{r}_2(\BB)$. Hence there are partitions
$\{a_{i,j} : i\in 2\}$ of $a_j$ which 2-reap
$\{a\wedge a_j : a\in A\}$ for each $j \in 2$. Consequently, $\{a_{i,j}
: \{i,j\}\subseteq 2\}$ is a 4-partition which $(4,4)$-reaps $A$.
Continuing in this manner one obtains an $i$-partition which $(i,i)$-reaps $A$
and hence, also $(i,j)$-reaps $A$.
 
Statement (5) is an immediate consequence of (4) and (3). In 
particular, (4) implies that $\frak{r}_{n,n}(\BB) \geq
\frak{r}_2(\BB)$ while  (3)
yields that  $\frak{r}_{n,m}(\BB) \leq
\frak{r}_2(\BB)$ for every Boolean algbera $\BB$ provided that
$n/(m-1)\geq 2$.
\stopproof

The question of whether it is possible to have a Boolean algebra with
three different reaping invariants was unresolved for some time.
The main theorem of this section, Theorem~\ref{suc},
 will show that this is not possible
and that even more restrictions apply. On the road to establishing
this it will prove to be useful to have the following definition of
{\em polarized partition relation}.
\begin{defin}
The collection of symbols
\[\left( \begin{array}{c} m \\ n \end{array} \right) \not\rightarrow
\left( \begin{array}{c} j \\ i \end{array} \right)^{1,1}_{k,q}\]
is defined to mean that there is a function, $h$, from
$n\times m$ into $k$ such that, for all $a \in [n]^i$ and $b\in [m]^j$
the function\label{ppp} $h\restricts a\times b$  takes on more than $q$ values.
\end{defin}
 
\begin{lemma}
\label{polpart}    For each prime integer $n$ and $k < n$
\[\left( \begin{array}{c} n \\ n \end{array} \right) \not\rightarrow
\left( \begin{array}{c} 2 \\ k \end{array} \right)^{1,1}_{n,k}\]
\end{lemma}
{\bf Proof: }
    Define the function $h : n\times n\rightarrow n$ by
$h(i,j) =  i + j \mod n$.\stopproof
 
\begin{lemma}
    Suppose that $i$, $j$, $k$, $m$, $n$ and $q$
are  integers such that \label{hlemm}
\[\left( \begin{array}{c} m \\ n \end{array} \right) \not\rightarrow
\left( \begin{array}{c} j \\ i \end{array} \right)^{1,1}_{k,q}\]
Then any
Boolean algebra $\Bbb{B}$ which satisfies that
$ \frak{r}_{n,i}(\Bbb{B}) > \frak{r}_{k,q+1}(\Bbb{B}) $ also satisfies
that $\frak{r}_{m,j}(\Bbb{B}) \leq \frak{r}_{k,q+1}(\Bbb{B})^+\ $.
 \end{lemma}
 
{\bf Proof: }
    Let $h: m\times n \rightarrow k$ witness that
\[\left( \begin{array}{c} m \\ n \end{array} \right) \not\rightarrow
\left( \begin{array}{c} j \\ i \end{array} \right)^{1,1}_{k,q}\]
Let $\Bbb{B}$ be a Boolean algebra and assume that $
\frak{r}_{n,i}(\Bbb{B}) \geq \kappa =  \frak{r}_{k,q+1}(\Bbb{B})^+$.
Choose a set $A \in [\Bbb{B}]^{\frak{r}_{k,q+1}(\Bbb{B})}$
which cannot be
$(k,q+1)$-reaped. Let $\Bbb{B}_0 $ be the algebra generated by $A$.
Recursively choose $n$-partitions $\{b(\alpha,\ell) : \ell < n\}$
which $(n,i)$-reap $\Bbb{B}_{\alpha}$ where, for $\beta < \kappa$,
 $\Bbb{B}_{\beta+1} $ is the
algebra generated by $\Bbb{B}_{\beta} \cup \{ b(\beta,\ell) : \ell < n\}$
and if $\beta $ is a limit ordinal, then $\Bbb{B}_\beta =
\bigcup_{\gamma < \beta}\Bbb{B}_\gamma$.
Finally, assuming that $\frak{r}_{m,j}(\Bbb{B}) > \kappa$, choose an
 $m$-partition,
$\{b(\kappa,\ell) : \ell < m\}$, which $(m,j)$-reaps $\Bbb{B}_\kappa$.
 
For each $\alpha < \kappa$ and $\ell < k$ define
\[  c(\alpha,\ell) = \bigvee \left\{ b(\alpha,\xi) \wedge
b(\kappa,\zeta) : h(\xi,\zeta) = \ell\right\}\ .\]
Since $\{ c(\alpha,\ell) : \ell < k\}$ is a $k$-partition of $\Bbb{B}$,
there are $\{\ell^e_\alpha: e\in q\} \subseteq k$ and
an element $a_\alpha \in A$ such
that $a_\alpha \leq \vee\{c(\alpha,\ell^e_\alpha) : e\in q\}$.
Now there are $\alpha_1
< \alpha_2 < \ldots < \alpha_i < \kappa$ so that $\ell^e_{\alpha_1} =
\ell^e_{\alpha_2} = \ldots = \ell^e_{\alpha_i} = \ell^e$, for each
$e\in q$, and $a_{\alpha_1} =
a_{\alpha_2} = \ldots = a_{\alpha_i} = a$.
 
Now choose inductively $I_{\eta}
< n$ so that
\begin{itemize}
    \item $I_{\eta} \neq I_{\eta'}$ if $\eta \neq \eta'$
\item $a \wedge b(\alpha_1,I_1)\wedge b(\alpha_2,I_2)\wedge\ldots\wedge
b(\alpha_i,I_i) > 0$
\end{itemize}
This is easily done because if $I_1, I_2, \dots I_\eta$ have been
chosen so that $\eta < i$ and
$a \wedge b(\alpha_1,I_1)\wedge b(\alpha_2,I_2)\wedge\ldots\wedge
b(\alpha_{\eta},I_{\eta}) > 0$ then there is some $w\in [n]^i$ such that
$a \wedge b(\alpha_1,I_1)\wedge b(\alpha_2,I_2)\wedge\ldots\wedge
b(\alpha_{\eta},I_{\eta})\wedge b(\alpha_{\eta+1},y) > 0$ for each $y\in w$
by virtue of the fact that $\Bbb{B}_{\alpha_{\eta+1}}$ is $(n,i)$-reaped by
$\{b(\alpha_{\eta+1},t) : t \in n \}$. It follows that it is
possible to choose $I_{\eta+1}\in w\setminus\{I_1, I_2, \ldots, I_\eta\}$.
 
  Since $a \wedge b(\alpha_1,I_1)\wedge b(\alpha_2,I_2)\wedge\ldots\wedge
b(\alpha_i,I_i) \in \Bbb{B}_{\kappa}$
it follows that  there is $u\in [m]^j$ so that
$a \wedge b(\alpha_1,I_1)\wedge b(\alpha_2,I_2)\wedge\ldots\wedge
b(\alpha_i,I_i)\wedge b(\kappa,y) > 0$ for each $y \in u$
by virtue of the fact that $\Bbb{B}_{\kappa}$ is $(m,j)$-reaped by
$\{b(\kappa,y) : y \in m \}$.
 
It will now be shown that $h(x,y) \in \{\ell^e : e\in q\}$
for any $x\in \{I_1, I_2, \dots ,I_i\}$ and
$y\in u$. To see this let $x = I_t$ and notice that $$a\wedge
b(\alpha_t,I_t)\wedge b(\kappa,y) \geq
a \wedge b(\alpha_1,I_1)\wedge b(\alpha_2,I_2)\wedge\ldots\wedge
b(\alpha_j,I_j)\wedge b(\kappa, y) > 0$$ Since $a \leq
\vee\{c(\alpha_t,\ell^e) : e\in q\}$ it follows that
$0 < a \wedge b(\alpha_t,I_t)\wedge b(\kappa,y)
< \vee\{c(\alpha_t, \ell^e) : e\in q\}$ and hence that
$h(I_t,y) = \ell^e$ for some $e\in q$.
This is a contradiction to the hypothesis on $h$
that its range on any $i \times j$ rectangle has more than $q$ points
in it.
\stopproof
\begin{corol}
  If $2\leq i \leq n\leq k$ then
$\frak{r}_{n,i}(\Bbb{B}) \leq \frak{r}_k
(\Bbb{B})^+$.
\end{corol}
\proof If $n = i$ then this follows from (4) and (5) of 
Lemma~\ref{mono} --- hence it may be assumed that $i < n$. Therefore  it
follows from Lemma~\ref{polpart} that
\[\left( \begin{array}{c} n \\ n \end{array} \right) \not\rightarrow
\left( \begin{array}{c} 2 \\ i \end{array} \right)^{1,1}_{n,i}\]
and consequently
\[\left( \begin{array}{c} n \\ n \end{array} \right) \not\rightarrow
\left( \begin{array}{c} i \\ i \end{array} \right)^{1,1}_{k,1}\]
because of simple monotonicity properties of polarized partition relations.
Finally let $n= m$ and $i= j$ and $ q= 1$ in Lemma~\ref{hlemm}. \stopproof

\begin{theor}
   The inequality $\frak{r}_{n}(\Bbb{B}) 
\leq \frak{r}_2(\Bbb{B})^+ $ \label{suc} holds for every integer $n$ and 
Boolean algebra $\Bbb{B}$.
\end{theor}
{\bf Proof} From Lemma~\ref{mono} (2) it suffices to show this only
in the case when $n$ is a prime. Let $\BB$ be a Boolean algebra.
If $\frak{r}_{n,2}(\BB) = \frak{r}_2(\BB)$ there is nothing to prove so it may
be assumed that there is some $k$ such that $\frak{r}_{n,k}(\BB) \neq
\frak{r}_2(\BB)$ --- let $k$ be the greatest such integer.
>From Lemma~\ref{mono} (5),
it follows that $k < n$ and from (4) of that same lemma it
follows that $\frak{r}_{n,k}(\BB) > \frak{r}_2(\BB)$.
>From Lemma~\ref{polpart} it follows that
\[\left( \begin{array}{c} n \\ n \end{array} \right) \not\rightarrow
\left( \begin{array}{c} 2 \\ k \end{array} \right)^{1,1}_{n,k}\]
and hence, from Lemma~\ref{hlemm} it  may be concluded that
$\frak{r}_{n,2}(\BB) \leq \frak{r}_{n,k+1}(\BB)^+$ because
the other hypothesis of Lemma~\ref{hlemm} --- namely that
$\frak{r}_{n,k}(\BB) > \frak{r}_{n,k+1}(\BB)$ ---
follows from the maximality assumption on
$k$. But, again because of the maximality of $k$, 
$\frak{r}_2(\BB) = \frak{r}_{n,k+1}(\BB)$ and hence
$\frak{r}_{n,2}(\BB) \leq \frak{r}_{2}(\BB)^+$.
\stopproof

The following immediate corollary to Theorem~\ref{suc}, together with
Proposition \ref{balc:1}, shows that the minimum $\pi$-character of a
point in the Stone space of any Boolean algebra $\Bbb{B}$ is bounded
above by $\frak{r}_2(\Bbb{B})^+$.
\begin{corol}
    For every Boolean algebra $\Bbb{B}$, \label{exp}
$\frak{r}_\omega(\Bbb{B}) \leq {\frak{r}_2(\Bbb{B})}^+$.
\end{corol}
\section{Finite Combinatorics}
It is the purpose of this section to provide some evidence for
the following conjecture:
For any pair of integers $i$ and $j$ there is an integer $n$ such that
 $\frak{r}_{n}(\BB) = \frak{r}_{i,j}(\BB)$ for any Boolean algebra $\BB$.
Indeed it is reasonable to conjecture that 
for any pair of integers $i \geq j \geq 2$  and any Boolean algebra $\BB$
$$\frak{r}_{n}(\BB) = \frak{r}_{i,j}(\BB)$$
where $n$ is the least integer greater than or equal to $i/(j-1)$.
Although this conjecture remains unproved it will be shown in this section
to be true when
$n = 2$ or $i \leq 8$ and in many other cases. Lemmas \ref{remath1} and \ref{remath2} will be used to do this.
\begin{lemma}
If $\BB$ is any Boolean algebra then \label{remath1}
$\frak{r}_{3,2}(\BB) = \frak{r}_{5,3}(\BB) = \frak{r}_{6,3}(\BB)$.
\end{lemma}
{\bf Proof: }
>From  Theorem~\ref{ijkl} it follows that in order
to show that $\frak{r}_{6,3}(\BB) \leq \frak{r}_{3,2}(\BB)$
it must be shown that $P(6,3,3,1)$ fails. This is easy since
the colouring of $\mu(T(1,6))$ which partitions $\mu(T(1,6))$ into three pairs
is a counterexample to $P(6,3,3,1)$. 

To show that  $\frak{r}_{3,2}(\BB) \leq \frak{r}_{5,3}(\BB)$
it must be shown that $P(3,5,2,2)$ fails.
Define $$\chi : \mu(T(3,3))\rightarrow 5$$ by
$\chi(\sigma) = \sum_{i\in 3}\sigma(i)\mod 5$ and check that this is
a counterexample.

That $\frak{r}_{5,3}(\BB) \leq \frak{r}_{6,3}(\BB)$ follows from 
Lemma~\ref{mono}.
\stopproof

\begin{lemma} If $\BB$ is any Boolean algebra \label{remath2}then
$\frak{r}_{3,2}(\BB) = \frak{r}_{7,4}(\BB) = \frak{r}_{8,4}(\BB)
= \frak{r}_{9,4}(\BB)$.
\end{lemma}
{\bf Proof: }
As in Lemma~\ref{remath1},
to show that $\frak{r}_{9,4}(\BB) \leq \frak{r}_{3,2}$
it must be shown that $P(9,3,4,1)$ fails. This is easy since
the colouring of $\mu(T(1,9))$ which partitions $\mu(T(1,9))$ 
into three triples
is a counterexample to $P(9,3,4,1)$. 

To show that  $\frak{r}_{3,2}(\BB) \leq \frak{r}_{9,4}(\BB)$
it must be shown that $P(3,7,2,3)$ fails.
Define $$\chi : \mu(T(3,3))\rightarrow 7$$ by
$\chi(\sigma) = \sum_{i\in 3}\sigma(i)\mod 7$ and check that this is
a counterexample.

That $\frak{r}(\BB)_{7,4} \leq \frak{r}(\BB)_{8,4}\leq \frak{r}_{9,4}(\BB)$ 
follows from Lemma~\ref{mono}.
\stopproof

The truth of the following conjecture, which was mentioned at the
beginning of this section can now be established in the cases $n= 2$
or $i\leq 8$:{\em For any pair of integers $i \geq j \geq 2$  
and any Boolean algebra $\BB$
$$\frak{r}_{n}(\BB) = \frak{r}_{i,j}(\BB)$$
where $n$ is the least integer greater than or equal to $i/(j-1)$.}\/
 First, if $n=2$ this is a direct consequence
of Lemma~\ref{mono} (5). If $j=2$ then the conjecture is a consequence of the 
definition of $\frak{r}_{i,2}(\BB)$. So, if $i\leq 8$ then from 
Lemma~\ref{remath2} it follows that $$\frak{r}_{8,4}(\BB) =
\frak{r}_{7,4}(\BB) = \frak{r}_{3,2}(\BB) = \frak{r}_{3}(\BB)$$
while from
Lemma~\ref{remath1} it follows that $$\frak{r}_{5,3}(\BB) =
\frak{r}_{6,3}(\BB) = \frak{r}_{3,2}(\BB) = \frak{r}_{3}(\BB)$$
so only $\frak{r}_{8,3}(\BB)$ and $\frak{r}_{7,3}(\BB)$ need be considered.
The following result of C. Laflamme \cite{lafl} deals with this case.
\begin{theor}
If $k \geq 2m-1$ then $P(m,k,2,2)$ f\label{laflam}ails and hence $\frak{r}_{m,2}(\BB)\leq
\frak{r}_{k,3}(\BB)$.
\end{theor}

The point is that for any Boolean algebra $\BB$ the reaping numbers $\frak{r}_{
7,3}(\BB)$ and $\frak{r}_{8,3}(\BB)$ are both less than or equal to 
$\frak{r}_{4,2}(\BB)$ by Lemma~\ref{mono} (5). The opposite inequalities are
consequences of Theorem~\ref{laflam}. Hence  $\frak{r}_{8,3}(\BB) =
\frak{r}_{4}(\BB) =  \frak{r}_{7,3}(\BB)$. This completes the first seven row of the following table in which the entry $n$ in row $i$ and column $j$ 
signifies that  $\frak{r}_{i,j}(\BB) =  \frak{r}_{n}(\BB)$ for every Boolean algebra $\BB$.

\bigskip

\begin{center}
    \begin{tabular}{|c|c|c|c|c|c|c|c|c|}
\hline
  & 2 & 3 & 4 & 5 & 6 & 7 & 8 & 9 \\ \hline
2 & 2 &   &   &   &   &   &   &   \\ \hline
3 & 3 & 2 &   &   &   &   &   &   \\ \hline
4 & 4 & 2 & 2 &   &   &   &   &   \\ \hline
5 & 5 & 3 & 2 & 2 &   &   &   &   \\ \hline
6 & 6 & 3 & 2 & 2 & 2 &   &   &   \\ \hline
7 & 7 & 4 & 3 & 2 & 2 & 2 &   &   \\ \hline
8 & 8 & 4 & 3 & 2 & 2 & 2 & 2 &   \\ \hline
9 & 9 & 5 & ? & 3 & 2 & 2 & 2 & 2 \\ \hline
     \end{tabular}\end{center}

\medskip

By using techniques similar in spirit to those of Lemmas \ref{remath1}
and \ref{remath2} it is possible to prove
the following.
\begin{lemma}If $\BB$ is an Boolean algebra then
\label{remath3}$\frak{r}_{3,2}(\BB) = \frak{r}_{9,5}(\BB) = \frak{r}_{12,5}(\BB)$.
\end{lemma}

This, together with Theorem~\ref{laflam} and Lemma~\ref{mono}, allow all but one of the entries in the last row of the table to be filled in.
The question mark indicates an open problem.
\begin{quest}
Does there exist a Boolean algebra $\BB$ such that
$\frak{r}_{9,4}(\BB)\neq\frak{r}_{3}(\BB)$?
\end{quest}
It is also possible to establish a simple arithmetic condition
which is equivalent to  $P(m,k,n,1)$.
\begin{lemma}
The property $P(m,k,n,1)$ fails  if and only\label{asdf} if $m < kn - k + 1$.
\end{lemma}
{\bf Proof: }
Colour $\mu(T(1,m))$ with $k$ colours in such a way that each colour gets used
at most $n-1$ times. It is possible to do this because
 $k(n-1) \geq m$. Any $n$ branching subtree gets at least
2 colours.

Suppose $kn - k + 1 \leq m$. It will be shown by induction on $h\in\omega$
that for  any colouring $\chi : \mu(T(h,m))
\rightarrow k$ there is a monochromatic $n$-branching subtree; in other
words, an $n$-branching subtree $S\subseteq T(h,m)$ such that $\mu(S)\subseteq
\mu(T(h,m))$ and $\chi\restriction\mu(S)$ is constant.
If $h=1 $ a simple pigeonhole argument can be applied.
Otherwise $h = h'+1$ and the induction hypothesis can be used
to find $c(i)\in k$ and an $n$-branching subtree $S_i\subseteq
T(h\setminus 1,m)$ such that $\mu(S_i) = \mu(T(h\setminus 1,m))$ and
$\chi(\sigma) = c(i)$ for each $\sigma\in \mu(S_i)$.
Again use a pigeonhole argument to find $X\in [m]^n$ and $c\in k$
such that $c(i) = c$ for $i \in X$ and then let 
$S= \{\sigma : \sigma(0)\in X\AND \sigma\restriction h\setminus
1\in S_{\sigma(0)}\}$.\stopproof

\begin{corol}
For any integer $k\geq 2$ 
there is some Boolean algebra $\BB$ 
such that $\frak{r}_{k+1,2}(\BB) \not\leq \frak{r}_{k,2}(\BB)$.
\end{corol}
{\bf Proof: } Let $n = 2$ in Lemma~\ref{asdf}.\stopproof

\begin{corol}
For any integer $k\geq 2$ 
there is some Boolean algebra $\BB$ 
such that $\frak{r}_{2k+1,3}(\BB) \not\leq \frak{r}_{k,2}(\BB)$.
\end{corol}
{\bf Proof: } Let $n = 3$ in Lemma~\ref{asdf}.\stopproof

\begin{corol} There are Boolean algebras $\BB_1$, $\BB_2$ and 
\label{qtp}$\BB_3$ such that
\begin{itemize}
\item  $\frak{r}_{9,5}(\BB) \not\leq \frak{r}_{2,2}(\BB)$.
\item  $\frak{r}_{10,4}(\BB) \not\leq \frak{r}_{3,2}(\BB)$.
\item  $\frak{r}_{9,3}(\BB) \not\leq \frak{r}_{4,2}(\BB)$.
\end{itemize}
\end{corol}
{\bf Proof: } Use Lemma~\ref{asdf}.\stopproof

It is possible to introduce a partial order $\prec$ on
$\omega\time\omega$ by defining $(n,m)\prec (i,j)$ if and only if
 $\frak{r}_{n,m}(\BB) \leq \frak{r}_{i,j}(\BB)$ for every  Boolean
algebra $\BB$.  Obviously $\prec$ is a transitive relation and
so it induces a partial order on equivalence classes. 
As an example, Lemma~\ref{remath1} shows that
$(3,2)$, $(5,3)$ and $(6,3)$ are all $\prec$ equivalent,
 Lemma~\ref{mono} shows that $(2,2)$ is $\prec$ minimal and
Corollary~\ref{qtp} shows that $(9,3)\not\prec (4,2)$. There are
various open questions about the partial order induced by $\prec$.
\begin{quest}
Is the $\prec$ interval between $(2,2)$ and $(3,2)$ empty?
\end{quest}
\begin{quest}
    What is the order type of $(\omega\times\omega,\prec)$?
\end{quest}
It is not even known if $(\omega\times\omega,\prec)$ is dense,
well founded or linear.
The simplest open question about the property $P(i,j,k,m)$ can be phrased as follows.
\begin{quest}
Does $(11,6) \prec (3,2)$ hold?
\end{quest}
\makeatletter \renewcommand{\@biblabel}[1]{\hfill#1.}\makeatother
\renewcommand{\bysame}{\leavevmode\hbox to3em{\hrulefill}\,}

\end{document}